\documentclass[12pt]{article} 
\usepackage{amssymb,amsmath}
\usepackage{amsthm}
\usepackage{bm}%
\usepackage{stmaryrd}
\usepackage {arydshln}
\usepackage{xcolor}
\usepackage{color}
\usepackage{comment}

  \def\cc{{\cal C}}
\def\cd{{\cal D}}  \def\cf{{\cal F}}
\def\cg{{\cal G}} \def\ch{{\cal H}} \def\ci{{\cal I}}
 \def\ck{{\cal K}} 
  
   \def\cs{{\cal S}}

\def\tcd{\tilde{\cal D}}

\newcommand{\eqna}[1]{\begin{subequations} \label{#1}
\begin{eqnarray}}
\def\eena{\end{eqnarray}
\end{subequations}}

 \def\half{{\textstyle{1\over2}}}
\def\a{\alpha}

\def\llra{\longleftrightarrow}

\def\beq{\begin{equation}}
\def\eqn#1{\beq\label{#1}}
\def\eeq{\end{equation}}
\def\ee{\end{equation}}
\def\lra{\longrightarrow}
\def\cg{{\cal G}}\def\ci{{\cal I}}
\def\nn{\nonumber}

\def\ket#1{\left| #1\right\rangle}
\def\hf{\frac{1}{2}}
\def\Gau#1{\left\lfloor #1 \right\rfloor}

\def\idos{intertwining differential operators}

\def\nn{\nonumber}
\def\nt{\noindent}

\def\bbc{\mathbb{C}}

\def\bbr{\mathbb{R}}
\begin{document}

\begin{center}
{\large\bf
  Invariant differential operators for the Jacobi algebra $ {\cal G}_2$
}

\vspace{10mm}

{{\large\bf  N. Aizawa$^1$, ~V.K. Dobrev$^2$}
	\\[10pt]
$^1$Department of Physical Science, Osaka Prefecture University, \\
Nakamozu Campus, Sakai, Osaka 599-8531, Japan
\\[10pt]
$^2$Institute of Nuclear Research and Nuclear Energy,\\
 Bulgarian Academy of Sciences, \\
72 Tsarigradsko Chaussee, 1784 Sofia, Bulgaria}

\end{center}

\vspace{10mm}

\begin{abstract}
In the present paper we construct explicitly the \idos\ for the     Jacobi algebra
$\cg_2$.  For the construction we use the singular vectors of the Verma modules over $\cg_2$ which we have constructued earlier. We construct the function spaces on which the operators act. We {\color{black} find} two versions of the
left (representation) action and the right action. {\color{black} These actions are combined with} the singular vectors to provide the \idos.
 \end{abstract}

\vspace{10mm}

\section{Introduction}

Consider a Lie group ~$G$, e.g., the Lorentz, Poincar\'e, conformal
groups, and differential equations
\eqn{inv0} \ci~f ~~=~~ j \ee
which are $G$-invariant. These
play a very important role in the description of physical
symmetries - recall, e.g., the early examples of Dirac, Maxwell, d'Allembert,
equations and nowadays the latest applications
of (super-)differential operators in conformal field theory,
supergravity, string theory, see e.g. \cite{Recent}. Naturally, it is important to construct systematically
such invariant equations and operators.

To recall the notions, consider a Lie group $G$ and
two representations $T,T'$
acting in the representation spaces $C,C'$, which may be Hilbert,
Fr\'echet, etc. An {\it invariant (or intertwining) operator} ~$\ci$~ for these two
representations is a continuous linear map
\eqn{ima0} \ci ~:~ C \lra C' \ee
such that
\eqn{int0}  T'(g) \circ \ci ~~=~~ \ci \circ T(g) ~, \quad \forall g\in
G ~. \ee
Then    we say that the equation \eqref{inv0} is a
$G$ - {\it invariant equation}.
Note that ~$\ker\ci$, ~im$\, \ci$~ are invariant subspaces of $C$,
$C'$, resp.

If $G$ is semisimple then there exist canonical ways
for the construction of the \idos, cf., e.g., \cite{Dob88,VKD1}.
 In this method there is a correspondence between invariant differential operators and singular vectors of Verma modules
over the (complexified) Lie algebra in consideration.

The procedure may be applied for more   general classes of Lie groups. For instance, it was applied to the
Schr\"odinger group \cite{Nie,Hag} in, e.g., \cite{DDM,AiDo}.

This is what we try to do in the present paper for the case of $\cg_2$.

\section{Preliminaries}
\setcounter{equation}{0}

The procedure that we shall follow requires first that we find the singular vectors of the Verma modules over  $\cg_2$.
This task was fulfilled in \cite{AiDoDo}. Furthermore there are given all necessary details, and thus, we can present the preliminaries in a shorter fashion.

The Jacobi algebra is the semi-direct sum $\cg_n:= \ch_n\niplus 
sp(n,\bbr)_{\bbc}$ \cite{EiZa,BeSc}. The Heisenberg algebra
$\ch_n$ is generated by the boson creation (respectively,
annihilation) operators~${a}_i^{+}$~(${a}^-_i$),~$i,j
=1,\dots,n$, which verify the canonical commutation relations
\eqn{heis} \big[a^-_i,a^{+}_j\big]=\delta_{ij}, \qquad [a^-_i,a^-_j]
= \big[a_i^{+},a_j^{+}\big]= 0 . \ee
$\ch_n$ is an
ideal in $\cg_n$, i.e., $[\ch_n,\cg_n]=\ch_n$,
determined by the commutation relations (following the notation of \cite{Berc}):
\eqna{haspn}
&&\big[a^{+}_k,K^+_{ij}\big] = [a^-_k,K^-_{ij}]=0, \\
&&{} [a^-_i,K^+_{kj}] = \tfrac{1}{2}\delta_{ik}a^{+}_j+\tfrac{1}{2}\delta_{ij}a^{+}_k ,\qquad
 \big[K^-_{kj},a^{+}_i\big] = \tfrac{1}{2}\delta_{ik}a^-_j+\tfrac{1}{2}\delta_{ij}a^-_k , \\
&& \big[K^0_{ij},a^{+}_k\big] = \tfrac{1}{2}\delta_{jk}a^{+}_i,\qquad
\big[a^-_k,K^0_{ij}\big]= \tfrac{1}{2}\delta_{ik}a^-_{j} .
\eena
 $K^{\pm,0}_{ij}$ are the generators of the $\cs_n ~\equiv~ sp(n,\bbr)_{\bbc}$ algebra:
\eqna{baspn}
&& [K_{ij}^-,K_{kl}^-] = [K_{ij}^+,K_{kl}^+]=0 , \qquad 2\big[K^-_{ij},K^0_{kl}\big] = K_{il}^-\delta_{kj}+K^-_{jl}\delta_{ki}\label{baza23}, \\
&& 2[K_{ij}^-,K_{kl}^+] = K^0_{kj}\delta_{li}+
K^0_{lj}\delta_{ki}+K^0_{ki}\delta_{lj}+K^0_{li}\delta_{kj}\\
&& 2\big[K^+_{ij},K^0_{kl}\big] = -K^+_{ik}\delta_{jl}-K^+_{jk}\delta_{li},\quad
 2\big[K^0_{ji},K^0_{kl}\big] = K^0_{jl}\delta_{ki}-K^0_{ki}\delta_{lj} . 
\eena

First, for simplicity, we introduce the following notations for the basis of ~$\cs_2$~:
\eqna{bas2} \cs^+ ~:~&&~ b^+_i ~\equiv~ K^+_{ii}\ , ~~i=1,2; \quad
c^+ ~\equiv~ K^+_{12}\ ,\quad d^+ ~\equiv K^0_{12} \\
\cs^- ~:~&&~ b^-_i ~\equiv~ K^-_{ii}\ , ~~i=1,2; \quad
c^- ~\equiv~ K^-_{12}\ ,\quad d^- ~\equiv K^0_{21} \\
\ck ~:~&&~ h_i ~\equiv~ K^0_{ii} \ , ~~i=1,2.
\eena

We need  also the triangular decomposition of $\cg_2$~:
\begin{align}
   {\cal G}_2^+ &:= \mathrm{l.s.}\{ \ a_i^+, \ b_i^+, \ c^+, \ d^+ \ \}, ~i=1,2,
   \nn \\
   {\cal G}_2^- &:= \mathrm{l.s.}\{ \ a_i^-, \ b_i^-, \ c^-, \ d^- \ \}, ~i=1,2,
   \nn \\
   {\cal K}_2 &:= \mathrm{l.s.}\{\ h_i, \ 1 \ \}, ~i=1,2.
\end{align}
{\color{black} Clearly,} {\color{black} the Abelian subalgebra $\ck$ is a Cartan subalgebra of $\cs_2$.
Furthermore,
~$\ck$~ plays the role  of
Cartan subalgebra for the whole algebra. Thus, we may treat the elements of  ${\cal G}_2^\pm$ as root subspaces w.r.t. $\ck$.
This may be explicated by the eigenvalues w.r.t. $(h_1,h_2)$ as follows:
\begin{align}  & a^\pm_1 ~:~ \pm(\half,0), ~~a^\pm_2 ~:~ \pm(0,\half), \\
& b^\pm_1 ~:~ \pm(1,0), ~~ b^\pm_2 ~:~ \pm(0,1), ~~  c^\pm ~:~ \pm (\half,\half), ~~ d^\pm ~:~ \pm (\half,-\half) \nn \end{align}
 We may also introduce the analogs of simple roots ~$\a_1,\a_2$~ which would correspond here to the
  generators ~$d^+,a^+_2$, resp. Then the correspondence  generators $\llra$ roots would be:
\eqn{rootgg}
(b^+_1, b^+_2, c^+, d^+ , a_1^+, a_2^+) ~\llra~ 
(2(\a_1+\a_2),2\a_2,\a_1+2\a_2 ,\a_1,\a_1+\a_2,\a_2) \eeq
}

{\color{black}
	We consider the lowest weight Verma modules over $\cg_2$ and found a complete list of singular vectors.
There are five types of singular vectors and they exist in Verma modules with a particular value of the lowest weight.
Explicit formula of them is found in \S 4.1 of \cite{AiDoDo}. 	
}

\bigskip

For the explicit construction of the \idos\ we need a parameter space. That would be some coset space of the
   Jacobi group $G$ as generated by ${\cal G}_2$.
Then we need its triangular decomposition $ G = G^+ K G^- $
and Borel subgroup $ B = K G^-.$

\bigskip
Now we can define the space of the right covariant functions:
\begin{equation}
   {\cal C}_{\Lambda} = \{ \ {\cal F} \in C^{\infty}(G) \ | \
     {\cal F}(g k g^-) = e^{\Lambda(H)}{\cal F}(g)\
   \}
\end{equation}
where $ g \in G, \ k = e^H \in K, \ g^- \in G^-, \  H \in {\cal K}$, $ \Lambda \in {\cal K}^*.$
{\color{black}
   Thus the functions of $ \cc_{\Lambda} $ are actually functions on $G/B$, or locally  on $G^+$.
}

\bigskip

Correspondingly we define the right action of ${\cal G}_2$ on ${\cal C}_{\Lambda}$:
\begin{equation}
    (\pi_R(X) {\cal F})(g) \doteq
    \left. \frac{d}{dt} {\cal F}(g \exp(tX)) \right|_{t=0},
    \qquad
    X \in {\cal G}_2, \ g \in G
\end{equation}
 and the left action of ${\cal G}_2$
\begin{equation}
    (\pi_L(X) {\cal F})(g) \doteq
    \left. \frac{d}{dt} {\cal F}(\exp(-tX)g) \right|_{t=0},
    \qquad
    X \in {\cal G}_2, \ g \in G
\end{equation}

In the next section we present these construction in explicit detail.

%
\section{Right and left actions on $G^+$} \label{SEC:RightAction}
\setcounter{equation}{0}

\subsection{Right action}

\bigskip
 For the elements $g^+$ of $ G^+$ we write:
\begin{equation}
   g^+ =  \exp\big(x_1 a_1^{\dagger} + x_2 a_2^{\dagger} \big)\,
   \exp\big(  y_1 b_1^+ + y_2 b_2^+ + z c^+ + w d^+    \big).
\end{equation}
It is important that there are only three non-vanishing relations among the generators of $ \cg_2^+:$
\begin{equation}
  [b_2^{\dagger}, d^{\dagger}] = -c^{\dagger},
  \qquad
  [a_2^{\dagger}, d^{\dagger}] = -\hf a_1^{\dagger},
  \qquad
  [c^{\dagger}, d^{\dagger}] = -\hf b_1^{\dagger}.  \label{NVcommR}
\end{equation}
Using these relations, it is easy to compute the right action of ${\cal G}_2^+:$
\begin{align}
  \pi_R(a_1^{\dagger}) &= \partial_{x_1},
  \nn \\
  \pi_R(a_2^{\dagger}) &= \partial_{x_2} + \frac{w}{2} \partial_{x_1},
  \nn \\
  \pi_R(b_1^{\dagger}) &= \partial_{y_1},
  \nn \\
  \pi_R(b_2^{\dagger}) &= \partial_{y_2} + \frac{w^2}{24} \partial_{y_1} + \frac{w}{2} \partial_z,
  \nn \\
  \pi_R(c^{\dagger}) &= \partial_z + \frac{w}{4} \partial_{y_1},
  \nn \\
  \pi_R(d^{\dagger}) &= \partial_w -\frac{1}{4} \left( z + \frac{y_2 w}{6} \right) \partial_{y_1} - \frac{y_2}{2} \partial_z.
  \label{RightAction}
\end{align}
As an example, we show the computation of $\pi_R(a_2^{\dagger})$.
First, noting \eqref{NVcommR} one may have
\begin{align}
	g^+ e^{ta_2^{\dagger}} &=
	e^A   \, e^{ta_2^{\dagger}} \, e^{-ta_2^{\dagger}} e^B  \, e^{ta_2^{\dagger}}
	\nonumber \\
	&= \exp(A+ta_2^{\dagger}) \exp\Big(B+\frac{tw}{2} a_1^{\dagger} \Big)
	\nonumber \\
	&= \exp\Big(A + t \big(a_2^{\dagger}+\frac{w}{2} a_1^{\dagger} \big) \Big) \exp(B)
\end{align}
where $  A:= x_1 a_1^{\dagger} + x_2 a_2^{\dagger}, B:= y_1 b_1^+ + y_2 b_2^+ + z c^+ + w d^+. $
It follows that
\begin{align}
	 (\pi_R(a_2^{\dagger}) {\cal F})(g^+) =
	\left. \frac{d}{dt} {\cal F}(g^+ e^{ta_2^{\dagger}}) \right|_{t=0}
	= \Big( \partial_{x_2} + \frac{w}{2} \partial_{x_1} \Big)  {\cal F}(g^+).
\end{align}

%
\subsection{Left action}

The left action is computed by using the Baker-Campbell-Hausdorff formula:
\begin{align}
  \ln e^X e^Y &= X + Y + \hf [X,Y] + \frac{1}{12} \big( (\mathrm{ad}X)^2 (Y) + (\mathrm{ad}Y)^2(X) \big)
  \nn \\
  & - \frac{1}{24} [Y, [X,[X,Y]]]
   - \frac{1}{720} \big( (\mathrm{ad}Y)^4(X) + (\mathrm{ad}X)^4(Y)  \big) + \cdots
\end{align}
where $ \mathrm{ad}X (Y) := [X,Y].$
We here present the final results and omit the computational details.

\bigskip\noindent
Left action of ${\cal G}_2^+:$
\begin{align}
  \pi_L(a_1^{\dagger}) &= -\partial_{x_1},
  \nn \\
  \pi_L(a_2^{\dagger}) &= -\partial_{x_2},
  \nn \\
  \pi_L(b_1^{\dagger}) &= -\partial_{y_1},
  \nn \\
  \pi_L(b_2^{\dagger}) &= -\partial_{y_2} - \frac{w^2}{24} \partial_{y_1} + \frac{w}{2} \partial_z,
  \nn \\
  \pi_L(c^{\dagger}) &= -\partial_{z} + \frac{w}{4} \partial_{y_1},
  \nn \\
  \pi_L(d^{\dagger}) &= -\partial_w - \frac{x_2}{2} \partial_{x_1}
  - \frac{1}{4} \left( z-\frac{y_2 w}{6} \right) \partial_{y_1}
  - \frac{y_2}{2} \partial_z.
\end{align}

\bigskip\noindent
Left action of ${\cal K}_2:$
\begin{align}
  \pi_L(h_1) &= -\frac{x_1}{2}\partial_{x_1} - y_1 \partial_{y_1} - \frac{z}{2} \partial_z - \frac{w}{2} \partial_w - \Lambda(h_1),
  \nn \\
  \pi_L(h_2) &= -\frac{x_2}{2}\partial_{x_2} - y_2 \partial_{y_2} - \frac{z}{2} \partial_z + \frac{w}{2} \partial_w - \Lambda(h_2),
  \nn \\
  \pi_L(1) &= -{\hat{\Lambda}}
  \end{align}
where $ {\hat{\Lambda}} $ is the  value of the central element of the Heisenberg algebra:
$ [a_i^-, a_j^+] = \delta_{ij} 1. $

\bigskip\noindent
Left action of ${\cal G}_2^-:$
\begin{align}
  \pi_L(a_1^-) &= -\left( y_1 + \frac{wz}{4} + \frac{y_2w^2}{12} \right) \partial_{x_1}
  - \frac{1}{2} \left( z + \frac{y_2w}{2} \right) \partial_{x_2} - x_1 {\hat{\Lambda}} ,
  \nn \\
  \pi_L(a_2^-) &= -\frac{1}{2} \left( z + \frac{y_2 w}{2} \right) \partial_{x_1}
  -y_2\, \partial_{x_2} - x_2 {\hat{\Lambda}} ,
  \nn \\
  \pi_L(b_1^-)
   &= x_1 \pi_L(a_1^-)
    - \left( y_1^2 + \frac{w^2}{96}\Big(z^2+y_2wz + \frac{y_2^2w^2}{12} \Big) \right) \partial_{y_1}
    -\frac{1}{4} \left( z + \frac{y_2w}{2} \right)^2 \partial_{y_2}
    \nn \\
    & -\left( y_1 z + \frac{w}{8} \Big( z^2 + \frac{2y_2 w z}{3} + \frac{y_2^2 w^2}{4} \Big) \right) \partial_{z}
    -w \left( y_1 -\frac{y_2w^2}{24}  \right)\partial_{w}
    \nn \\
    & + \frac{x_1^2}{2}{\hat{\Lambda}}   - 2\left( y_1 - \frac{y_2 w^2}{24} \right) \Lambda(h_1)
      - \frac{w}{2} \left( z + \frac{y_2 w}{2} \right) \Lambda(h_2),
   \nn \\
   \pi_L(b_2^-)
    &= x_2 \pi_L(a_2^-)  - \frac{y_2 w z}{12} \partial_{y_1}
    - y_2^2\, \partial_{y_2}
    \nn \\
    &-\frac{y_2}{2} \left( z + \frac{y_2 w}{2} \right) \partial_{z}
    -\left( z - \frac{y_2 w}{2} \right) \partial_{w}
    +\frac{x_2^2}{2} {\hat{\Lambda}}   - 2y_2\, \Lambda(h_2),
    \nn \\
    \pi_(c^-)
    &= \frac{x_2}{2} \pi_L(a_1^-) + \frac{x_1}{2} \pi_L(a_2^-)
    \nn \\
    & - \frac{1}{4} \left( y_1 \Big(z - \frac{y_2 w}{6} \Big) + \frac{y_2 w}{8} \Big( wz + \frac{y_2 w^2}{18} \Big) \right) \partial_{y_1}
    - \frac{y_2}{2} \left( z+ \frac{y_2w}{2}  \right) \partial_{y_2}
    \nn \\
    & -\hf \left( y_2 \Big( y_1 + \frac{wz}{4}+ \frac{5 y_2 w^2}{24} \Big) + \frac{z^2}{2}   \right) \partial_z
    - \left( y_1 + \frac{wz}{4} - \frac{y_2 w^2}{6} \right) \partial_w
    \nn \\
    & + \frac{x_1 x_2}{2}{\hat{\Lambda}}    - \hf \left( z - \frac{y_2 w}{2} \right) \Lambda(h_1)
      - \hf \left( z + \frac{3y_2 w}{2} \right) \Lambda(h_2),
    \nn \\
    \pi_L(d^-) &= -\frac{x_1}{2} \partial_{x_2} + \frac{w}{4} \left( y_1 - \frac{wz}{12} \right) \partial_{y_1}
    - \hf \left( z + \frac{y_2 w}{2} \right) \partial_{y_2}
    \nn \\
    &- \left( y_1 + \frac{y_2 w^2}{12} \right) \partial_z
     + \frac{w^2}{4} \partial_w + \frac{w}{2} \big( \Lambda(h_1) - \Lambda(h_2) \big).
\end{align}
It has been verified by direct computation (with MAPLE) that the left action given above is compatible with the defining commutation relations of $ {\cal G}_2.$


\section{Invariant differential operators: first version}
\setcounter{equation}{0}

First  we give the list of singular vectors that were found in \cite{AiDoDo}.
We denote the lowest weight vector of the Verma module by $ \ket{0} $ and the lowest weight by $\Lambda_k = \Lambda(h_k).$
The parameters $p^k $ and $q^3$ take a positive integer and the weight of the singular vector is denoted by $\Lambda'.$
\begin{enumerate}
 \renewcommand{\labelenumi}{(\roman{enumi})}
 \item $ \Lambda_1 - \Lambda_2 = \hf(1-p^1)$
 \begin{equation}
   \ket{v_s^{\Lambda'}} = (d^+)^{p^1} \ket{0}, \quad
   \Lambda' = \Lambda+p^1 (\delta_1-\delta_2). \label{SVcase1}
 \end{equation}
 \item $ \text{for all}\ \Lambda_1, \ \Lambda_2 = \frac{3}{4}- \frac{p^2}{2} $
 \begin{equation}
    \ket{v_s^{\Lambda'}} = (\hat{b}_2^+)^{p^2} \ket{0}, \quad
    \Lambda' = \Lambda + 2 p^2 \delta_2.
    \label{SVcase2}
 \end{equation}
 \item $ \Lambda_1 = \frac{5}{4}- \hf(p^3-q^3), \ \Lambda_2 = \frac{3}{4}-\hf q^3, \ (p^3 \neq q^3, p^3 \neq 2q^3)$
   \begin{enumerate}
     \item $ p^3 < q^3$
       \begin{equation}
          \ket{v_s^{\Lambda'}} =  \sum_{k=0}^{\Gau{p^3/2}} \sum_{n=0}^{p^3-2k} c(k,n)
          \ket{k, q^3-k-n, n, p^3-2k-n}
       \end{equation}
     \item $ q^3 < p^3 < 2q^3 $
       \begin{equation}
          \ket{v_s^{\Lambda'}} =
          \left(
            \sum_{k=0}^{p^3-q^3} \sum_{n=0}^{q^3-k} + \sum_{k=p^3-q^3+1}^{\Gau{p^3/2}} \sum_{n=0}^{p^3-2k}
          \right) c(k,n)
          \ket{k, q^3-k-n, n, p^3-2k-n}
       \end{equation}
     \item $ 2q^3 < p^3 $
       \begin{equation}
          \ket{v_s^{\Lambda'}} =  \sum_{k=0}^{q^3} \sum_{n=0}^{q^3-k} c(k,n)
          \ket{k, q^3-k-n, n, p^3-2k-n}
       \end{equation}
   \end{enumerate}
     where $\Lambda'$ and $c(k,n)$ are common for (a) (b) (c) and given by
     \begin{align}
        \Lambda' &= \Lambda + p^3 \delta_1 + (2q^3-p^3) \delta_2,
        \\
        c(k,n) &= \frac{ p^3! q^3! }{4^k k! n! (p^3-2k-n)! (q^3-k-n)!}.
     \end{align}
   \item $ \Lambda_1 + \Lambda_2 = 2 - \frac{p^4}{2}$
   \begin{align}
       \ket{v_s^{\Lambda'}} &=  \sum_{k=0}^{\Gau{p^4/2}} \sum_{n=0}^{p^4-2k}
       c(k,n) \ket{k,p^4-k-n,n,p^4-2k-n},
       \quad
       \Lambda' = \Lambda + p^4(\delta_1 + \delta_2),
       \\
       c(k,n) &= \frac{p^4!}{4^k k! n! (p^4-2k-n)!} \frac{\Gamma(2\Lambda_1+p^4-\frac{3}{2})}{\Gamma(2\Lambda_1+p^4-\frac{3}{2}-k-n)}.
   \end{align}
   \item $ \Lambda_1 = \frac{5}{4}- \frac{p^5}{2}, \ \text{for all}\ \Lambda_2 $
   \begin{align}
       \ket{v_s^{\Lambda'}} &=  \sum_{k=0}^{p^5} \sum_{n=0}^{p^5-k}
       c(k,n) \ket{k,p^5-k-n,n,2p^5-2k-n},
       \quad
       \Lambda' = \Lambda + 2p^5 \delta_1,
       \\
       c(k,n) &= \frac{(-1)^n}{4^k k! n!} \frac{p^5!}{(p^5-k-n)!}
      \frac{\Gamma(2\Lambda_2-p^5-\frac{3}{2}+2k+n)}{\Gamma(2\Lambda_2-p^5-\frac{3}{2})}.
   \end{align}
\end{enumerate}

Note that there is a change of basis w.r.t. \cite{AiDoDo}, namely:
\begin{align}
  \ket{k,\ell, n,m} :=  (\hat{b}_1^+)^k (\hat{b}_2^+)^{\ell} (\hat{c}^+)^n (d^+)^{m} \ket{0}
\end{align}
and
 \eqn{bchat}
   \hat{b}_k^+ := b_k^+ - \hf (a_k^+)^2, \qquad
   \hat{c}^+ := c^+ - \hf a_1^+ a_2^+.
\end{equation}

Then using the right action obtained in \S \ref{SEC:RightAction} we have  using \eqref{bchat}:
\eqn{rac} \begin{aligned}
  \pi_R(\hat{b}_1^+) &= \partial_{y_1} -\hf \partial_{x_1}^2,
  \cr
  \pi_R(\hat{b}_2^+) &=   \partial_{y_2} + \frac{w^2}{24} \partial_{y_1} + \frac{w}{2} \partial_z
        -\hf \Big(  \partial_{x_2} + \frac{w}{2} \partial_{x_1} \Big)^2,
   \cr
   \pi_R(\hat{c}^+) &=  \partial_z + \frac{w}{4} \partial_{y_1}
   -\hf \left( \partial_{x_1} \partial_{x_2} + \frac{w}{2} \partial_{x_1}^2 \right).
\end{aligned}
\ee

Thus the invariant differential operators are given by substituting the above right action in the singular vectors
given above:
\begin{enumerate}
 \renewcommand{\labelenumi}{(\roman{enumi})}
 \item $ \Lambda_1-\Lambda_2 = \hf(1-p^1)$
   \begin{equation} 
   \tcd_{(i)} =  \left(
       \partial_w -\frac{1}{4} \left( z + \frac{y_2 w}{6} \right) \partial_{y_1} - \frac{y_2}{2} \partial_z
     \right)^{p^1}
   \end{equation}
 \item $ \text{for all}\ \Lambda_1, \ \Lambda_2 = \frac{3}{4}- \frac{p^2}{2} $
   \begin{equation}
     \tcd_{(ii)} = \left(
        \partial_{y_2} + \frac{w^2}{24} \partial_{y_1} + \frac{w}{2} \partial_z
        -\hf
        \Big(  \partial_{x_2} + \frac{w}{2} \partial_{x_1} \Big)^2
      \right)^{p^2}
   \end{equation}
 \item $ \Lambda_1 = \frac{5}{4}- \hf(p^3-q^3), \ \Lambda_2 = \frac{3}{4}-\hf q^3, \ (p^3 \neq q^3, p^3 \neq 2q^3)$
   \begin{enumerate}
     \item $ p^3 < q^3$
       \begin{equation}
     \tcd_{(iii,a)} =     \sum_{k=0}^{\Gau{p^3/2}} \sum_{n=0}^{p^3-2k} c(k,n)\, {\cal P}(p^3,q^3,k,n)
       \end{equation}
     \item $ q^3 < p^3 < 2q^3 $
       \begin{equation}
        \tcd_{(iii,b)} =   \left(
            \sum_{k=0}^{p^3-q^3} \sum_{n=0}^{q^3-k} + \sum_{k=p^3-q^3+1}^{\Gau{p^3/2}} \sum_{n=0}^{p^3-2k}
          \right) c(k,n)\, {\cal P}(p^3,q^3,k,n)
       \end{equation}
     \item $ 2q^3 < p^3 $
       \begin{equation}
        \tcd_{(iii,c)} =   \sum_{k=0}^{q^3} \sum_{n=0}^{q^3-k} c(k,n)
          \, {\cal P}(p^3,q^3,k,n)
       \end{equation}
   \end{enumerate}
     where $ {\cal P}(p^3,q^3,k,n)$ and $c(k,n)$ are common for (a) (b) (c) and given by
     \begin{align}
        {\cal P}(p^3,q^3,k,n) &=
           \left( \partial_{y_1} -\hf \partial_{x_1}^2 \right)^k
           \left(
             \partial_{y_2} + \frac{w^2}{24} \partial_{y_1} + \frac{w}{2} \partial_z
             -\hf \Big(  \partial_{x_2} + \frac{w}{2} \partial_{x_1} \Big)^2
           \right)^{q^3-k-n}
           \nn\\
           & \times
           \left(
              \partial_z + \frac{w}{4} \partial_{y_1}
             -\hf \left( \partial_{x_1} \partial_{x_2} + \frac{w}{2} \partial_{x_1}^2 \right)
           \right)^n
           \nn \\
           & \times
           \left(
             \partial_w -\frac{1}{4} \left( z + \frac{y_2 w}{6} \right) \partial_{y_1} - \frac{y_2}{2} \partial_z
           \right)^{p^3-2k-n},
        \\
        c(k,n) &= \frac{ p^3! q^3! }{4^k k! n! (p^3-2k-n)! (q^3-k-n)!}.
     \end{align}
   \item $ \Lambda_1 + \Lambda_2 = 2 - \frac{p^4}{2}$
     \begin{align}
   \tcd_{(iv)} =     &  \sum_{k=0}^{\Gau{p^4/2}} \sum_{n=0}^{p^4-2k}
       \frac{p^4!}{4^k k! n! (p^4-2k-n)!} \frac{\Gamma(2\Lambda_1+p^4-\frac{3}{2})}{\Gamma(2\Lambda_1+p^4-\frac{3}{2}-k-n)}
       \nn \\
           & \qquad \times
           \left( \partial_{y_1} -\hf \partial_{x_1}^2 \right)^k
           \left(
             \partial_{y_2} + \frac{w^2}{24} \partial_{y_1} + \frac{w}{2} \partial_z
             -\hf \Big(  \partial_{x_2} + \frac{w}{2} \partial_{x_1} \Big)^2
           \right)^{p^4-k-n}
           \nn\\
           & \qquad  \times
           \left(
              \partial_z + \frac{w}{4} \partial_{y_1}
             -\hf \left( \partial_{x_1} \partial_{x_2} + \frac{w}{2} \partial_{x_1}^2 \right)
           \right)^n
           \nn \\
           &  \qquad \times
           \left(
             \partial_w -\frac{1}{4} \left( z + \frac{y_2 w}{6} \right) \partial_{y_1} - \frac{y_2}{2} \partial_z
           \right)^{p^4-2k-n}
     \end{align}
   \item $ \Lambda_1 = \frac{5}{4}- \frac{p^5}{2}, \ \text{for all}\ \Lambda_2 $
     \begin{align}
     \tcd_{(v)} =    &  \sum_{k=0}^{p^5} \sum_{n=0}^{p^5-k}
         \frac{(-1)^n}{4^k k! n!} \frac{p^5!}{(p^5-k-n)!}
      \frac{\Gamma(2\Lambda_2-p^5-\frac{3}{2}+2k+n)}{\Gamma(2\Lambda_2-p^5-\frac{3}{2})}
       \nn \\
           & \qquad \times
           \left( \partial_{y_1} -\hf \partial_{x_1}^2 \right)^k
           \left(
             \partial_{y_2} + \frac{w^2}{24} \partial_{y_1} + \frac{w}{2} \partial_z
             -\hf \Big(  \partial_{x_2} + \frac{w}{2} \partial_{x_1} \Big)^2
           \right)^{p^5-k-n}
           \nn\\
           & \qquad  \times
           \left(
              \partial_z + \frac{w}{4} \partial_{y_1}
             -\hf \left( \partial_{x_1} \partial_{x_2} + \frac{w}{2} \partial_{x_1}^2 \right)
           \right)^n
           \nn \\
           &  \qquad \times
           \left(
             \partial_w -\frac{1}{4} \left( z + \frac{y_2 w}{6} \right) \partial_{y_1} - \frac{y_2}{2} \partial_z
           \right)^{2p^5-2k-n}
    \end{align}

\end{enumerate}

%
\section{Final expressions for the \idos}
\label{FEidos}
\setcounter{equation}{0}

To simplify  our results we
make the following change of parameters:\\ $ (x_1,x_2,y_1,y_2,z,w) \to (\xi_1, \xi_2,\eta_1,\eta_2, \zeta, \omega)$,
namely:
\begin{alignat}{3}
  \xi_1 &= x_1 - \frac{w}{2}x_2, & \qquad
  \xi_2 &= x_2,
  \nn \\
  \eta_1 &= y_1 + \frac{w^2}{12} y_2 - \frac{w}{4}z, &
  \eta_2 &= y_2,
  & \qquad
  \zeta &= z - \frac{w}{2} y_2, \qquad \omega = w.
\end{alignat}
Inverse transform is given by
\begin{alignat}{3}
   x_1 &= \xi_1 + \frac{\omega}{2}\xi_2 ,
   & \qquad
   x_2 &= \xi_2,
   \nn \\
   y_1 &= \eta_1 + \frac{\zeta \omega}{4} + \frac{\eta_2 \omega^2}{24},
   & \qquad
   y_2 &= \eta_2,
   & \qquad
   z &= \zeta + \frac{\eta_2 \omega}{2}, \qquad w = \omega.
\end{alignat}
It follows that
\begin{align}
   \partial_{x_1} &= \partial_{\xi_1},
   \qquad
   \partial_{x_2} = \partial_{\xi_2} - \frac{\omega}{2} \partial_{\xi_1},
  \nn \\
   \partial_{y_1} &= \partial_{\eta_1},
    \qquad
   \partial_{y_2} = \partial_{\eta_2} + \frac{\omega^2}{12}\partial_{\eta_1} - \frac{\omega}{2}\partial_{\zeta},
   \qquad
   \partial_{z} = \partial_{\zeta} - \frac{\omega}{4}\partial_{\eta_1},
   \nn \\
   \partial_{w} &= \partial_{\omega} - \frac{1}{4}\left( \zeta - \frac{\eta_2\omega}{6} \right) \partial_{\eta_1} -  \frac{\eta_2}{2}\partial_{\zeta} - \frac{\xi_2}{2}\partial_{\xi_1}.
\end{align}
 Then,  the right action becomes:
\begin{align}
  \pi_R(a_k^+) &= \partial_{\xi_k},
  \qquad
  \pi_R(b_k^+) = \partial_{\eta_k},
  \qquad
  \pi_R(c^+) = \partial_{\zeta},
  \nn \\
  \pi_R(d^+) &= \partial_{\omega} - \frac{\zeta}{2}\partial_{\eta_1} - \eta_2 \partial_{\zeta} - \frac{\xi_2}{2}\partial_{\xi_1}.
\end{align}

The left action is also simplified and now reads as follows:
\begin{alignat}{2}
  \pi_L(a_1^+) &= -\partial_{\xi_1},
  & \qquad
  \pi_L(a_2^+) &= -\partial_{\xi_2} + \frac{\omega}{2}\partial_{\xi_1},
  \nn \\
  \pi_L(b_1^+) &= -\partial_{\eta_1},
  & \qquad
  \pi_L(b_2^+) &= -\partial_{\eta_2} - \frac{\omega^2}{4}\partial_{\eta_1} + \omega\partial_{\zeta},
  \nn \\
  \pi_L(c^+) &= -\partial_{\zeta} + \frac{\omega}{2}\partial_{\eta_1},
  &
  \pi_L(d^+) &= -\partial_{\omega},
\end{alignat}
\begin{align}
  \pi_L(h_1) &= -\frac{\xi_1}{2}\partial_{\xi_1} - \eta_1 \partial_{\eta_1} - \frac{\zeta}{2}\partial_{\zeta} - \frac{\omega}{2}\partial_{\omega} - \Lambda(h_1),
  \nn\\
  \pi_L(h_2) &= -\frac{\xi_2}{2}\partial_{\xi_2} - \eta_2 \partial_{\eta_2} - \frac{\zeta}{2}\partial_{\zeta} + \frac{\omega}{2}\partial_{\omega} - \Lambda(h_2),
  \nn \\
  \pi_L(1) &= -{\hat{\Lambda}}  ,
\end{align}
\begin{align}
   \pi_L(a_1^-) &= -\left( \eta_1 + \frac{\zeta\omega}{4}\right) \partial_{\xi_1}
   -\hf( \zeta + \eta_2 \omega) \partial_{\xi_2}
   -\left( \xi_1 + \frac{\xi_2 \omega}{2}\right) {\hat{\Lambda}} ,
   \nn \\
   \pi_L(a_2^-) &= - \frac{\zeta}{2}\partial_{\xi_1} - \eta_2 \partial_{\xi_2} - \xi_2 {\hat{\Lambda}} ,
   \nn \\
   \pi_L(b_1^-) &= \left(\xi_1 + \frac{\xi_2 \omega}{2}\right)  \pi_L(a_1^-)
     + \frac{\xi_2 \omega}{2}\left( \eta_1 + \frac{\zeta \omega}{4} \right) \partial_{\xi1}
     - \left(\eta_1^2-\frac{\zeta^2 \omega^2}{16}\right) \partial_{\eta_1}
   \nn\\
   &- \frac{1}{4}(\zeta + \eta_2 \omega)^2 \partial_{\eta_2}
    -  \left(\eta_1 + \frac{\zeta \omega}{4}\right)( \zeta\partial_{\zeta} + \omega\partial_{\omega}-2 \Lambda(h_1))
     \nn\\
    & + \hf \left( \xi_1 + \frac{\xi_2 \omega}{2}\right)^2 {\hat{\Lambda}}
      - \frac{\omega}{2}(\zeta + \eta_2 \omega) \Lambda(h_2),
   \nn \\
   \pi_L(b_2^-) &= \xi_2 \pi_L(a_2^-) + \frac{\xi_2 \zeta}{2} \partial_{\xi_1}
      + \frac{\zeta^2}{4}\partial_{\eta_1} - \eta_2^2 \partial_{\eta_2}
      - \zeta \partial_{\omega}
      + \frac{\xi_2^2}{2} {\hat{\Lambda}}-2\eta_2 \Lambda(h_2),
   \nn \\
   \pi_L(c^-) &= \frac{\xi_2}{2} \pi_L(a_1^-) + \hf \left( \xi_1 + \frac{\xi_2 \omega}{2} \right) \pi_L(a_2^-)
      +  \left( \eta_1 + \frac{\zeta \omega}{2} \right) \left( \frac{\xi_2}{2} \partial_{\xi_1} - \partial_{\omega} \right)
   \nn \\
   & + \frac{\zeta^2 \omega}{8} \partial_{\eta_1}
     - \hf (\zeta + \eta_2 \omega) ( \eta_2 \partial_{\eta_2} + \Lambda(h_2))
     - \frac{\zeta^2}{4}\partial_{\zeta}
   \nn \\
   & + \frac{\xi_2}{2}\left( \xi_1 + \frac{\xi_2 \omega}{2} \right) {\hat{\Lambda}}
     - \frac{\zeta}{2} \Lambda(h_1),
   \nn \\
   \pi_L(d^-) &= \frac{\omega \xi_1}{4}\partial_{\xi_1} - \hf \left(\xi_1+ \frac{\omega \xi_2}{2}\right) \partial_{\xi_2}
   + \frac{\omega \eta_1}{2}\partial_{\eta_1} - \hf (\zeta + \omega\eta_2)\partial_{\eta_2}
   \nn \\
   & - \eta_1 \partial_{\zeta} + \frac{\omega^2}{4}\partial_{\omega}
     + \frac{\omega}{2}( \Lambda(h_1) - \Lambda(h_2)).
\end{align}
$ \pi_L(b_k^-), \pi_L(c^-), \pi_L(d^-) $ have the following simpler expressions:
\begin{align}
   \pi_L(b_1^-) &= \xi_1 \pi_L(a_1^-) - \frac{\xi_1 \omega}{2} \pi_L(a_2^-)
      + \omega \pi_L(c^-) - \frac{\omega^2}{4} \pi_L(b_2^-)
   \nn \\
   & -\eta_1^2 \partial_{\eta_1} - \frac{\zeta^2}{4}\partial_{\eta_2} - \eta_1 \zeta \partial_{\zeta}
   + \frac{\xi_1^2}{2} {\hat{\Lambda}} -2 \eta_1 \Lambda(h_1),
   \nn\\
   \pi_L(b_2^-) &= -\xi_2 \eta_2 \partial_{\xi_2} + \frac{\zeta^2}{4}\partial_{\eta_1} - \eta_2^2 \partial_{\eta_2}
      - \zeta \partial_{\omega}
      - \frac{\xi_2^2}{2} {\hat{\Lambda}}-2\eta_2 \Lambda(h_2),
   \nn \\
   \pi_L(c^-) &= \frac{\xi_2}{2} \pi_L(a_1^-)
      + \hf \left( \xi_1 - \frac{\xi_2 \omega}{2}\right) \pi_L(a_2^-)
      + \frac{\omega}{2} \pi_L(b_2^-)
   \nn \\
   & + \frac{\xi_2 \eta_1}{2} \partial_{\xi_1} - \frac{\zeta \eta_2}{2} \partial_{\eta_2}
     - \frac{\zeta^2}{4} \partial_{\zeta} - \eta_1 \partial_{\omega}
       + \frac{\xi_1 \xi_2}{2} {\hat{\Lambda}}
       - \frac{\zeta}{2} (\Lambda(h_1) + \Lambda(h_2)),
   \nn \\
   \pi_L(d^-) &= - \frac{\omega}{2}\big( \pi_L(h_1)- \pi_L(h_2) \big)
     - \frac{\xi_1}{2}\partial_{\xi_2} - \frac{\zeta}{2}\partial_{\eta_2}
      -\eta_1 \partial_{\zeta} - \frac{\omega^2}{4}\partial_{\omega}.
\end{align}


 \bigskip

Finally, we pass to the "hat" basis:
\begin{align}
  \pi_R(\hat{b}_k^+) &= \partial_{\eta_k} -\hf \partial_{\xi_k}^2,
  \qquad
   \pi_R(\hat{c}^+) =  \partial_{\zeta}
   -\hf  \partial_{\xi_1} \partial_{\xi_2}.
\end{align}

Then the {\bf final} expressions for the invariant differential operators are:
\begin{enumerate}
 \renewcommand{\labelenumi}{(\roman{enumi})}
 \item $ \Lambda_1-\Lambda_2 = \hf(1-p^1)$
   \begin{equation}
    \cd_{(i)} ~=~  \left(
       \partial_{\omega} - \frac{\xi_2}{2} \partial_{\xi_1} - \frac{\zeta}{2}\partial_{\eta_1} - \eta_2 \partial_{\zeta}
     \right)^{p^1}
   \end{equation}
 \item $ \text{for all}\ \Lambda_1, \ \Lambda_2 = \frac{3}{4}- \frac{p^2}{2} $
   \begin{equation}
     \cd_{(ii)} ~=~  \left(
        \partial_{\eta_2} -\hf \partial_{\xi_2}^2
      \right)^{p^2}
   \end{equation}
 \item $ \Lambda_1 = \frac{5}{4}- \hf(p^3-q^3), \ \Lambda_2 = \frac{3}{4}-\hf q^3, \ (p^3 \neq q^3, p^3 \neq 2q^3)$
   \begin{enumerate}
     \item $ p^3 < q^3$
       \begin{align}
     \cd_{(iii,a)} ~=~      \sum_{k=0}^{\Gau{p^3/2}} \sum_{n=0}^{p^3-2k} &
           \frac{ p^3! q^3! }{4^k k! n! (p^3-2k-n)! (q^3-k-n)!}
           \left( \partial_{\eta_1} -\hf \partial_{\xi_1}^2 \right)^k
           \nn \\
           & \times
           \left(
             \partial_{\eta_2} -\hf \partial_{\xi_2}^2
           \right)^{q^3-k-n}
           \left(
              \partial_{\zeta}  -\hf  \partial_{\xi_1} \partial_{\xi_2}
           \right)^n
           \nn \\
           & \times
           \left(
               \partial_{\omega} - \frac{\xi_2}{2} \partial_{\xi_1}
               - \frac{\zeta}{2}\partial_{\eta_1} - \eta_2 \partial_{\zeta}
           \right)^{p^3-2k-n}
       \end{align}
     \item $ q^3 < p^3 < 2q^3 $
       \begin{equation}
        \cd_{(iii,b)} ~=~   \left(
            \sum_{k=0}^{p^3-q^3} \sum_{n=0}^{q^3-k} + \sum_{k=p^3-q^3+1}^{\Gau{p^3/2}} \sum_{n=0}^{p^3-2k}
          \right) \, {\cal P}(p^3,q^3,k,n)
       \end{equation}
     \item $ 2q^3 < p^3 $
       \begin{equation}
        \cd_{(iii,c)} ~=~   \sum_{k=0}^{q^3} \sum_{n=0}^{q^3-k} \, {\cal P}(p^3,q^3,k,n)
       \end{equation}
   \end{enumerate}
     where the summand $ {\cal P}(p^3,q^3,k,n)$  for (b) (c) is same as (a).
   \item $ \Lambda_1 + \Lambda_2 = 2 - \frac{p^4}{2}$
     \begin{align}
  \cd_{(iv)} ~=~      &  \sum_{k=0}^{\Gau{p^4/2}} \sum_{n=0}^{p^4-2k}
       \frac{p^4!}{4^k k! n! (p^4-2k-n)!} \frac{\Gamma(2\Lambda_1+p^4-\frac{3}{2})}{\Gamma(2\Lambda_1+p^4-\frac{3}{2}-k-n)}
       \nn \\
           & \qquad \times
           \left( \partial_{\eta_1} -\hf \partial_{\xi_1}^2 \right)^k
           \left(
             \partial_{\eta_2} -\hf \partial_{\xi_2}^2
           \right)^{p^4-k-n}
           \left(
              \partial_{\zeta}  -\hf  \partial_{\xi_1} \partial_{\xi_2}
           \right)^n
           \nn \\
           &  \qquad \times
           \left(
               \partial_{\omega} - \frac{\xi_2}{2} \partial_{\xi_1}
               - \frac{\zeta}{2}\partial_{\eta_1} - \eta_2 \partial_{\zeta}
           \right)^{p^4-2k-n}
     \end{align}
   \item $ \Lambda_1 = \frac{5}{4}- \frac{p^5}{2}, \ \text{for all}\ \Lambda_2 $
     \begin{align}
    \cd_{(v)} ~=~    &  \sum_{k=0}^{p^5} \sum_{n=0}^{p^5-k}
         \frac{(-1)^n}{4^k k! n!} \frac{p^5!}{(p^5-k-n)!}
      \frac{\Gamma(2\Lambda_2-p^5-\frac{3}{2}+2k+n)}{\Gamma(2\Lambda_2-p^5-\frac{3}{2})}
       \nn \\
           & \qquad \times
           \left( \partial_{\eta_1} -\hf \partial_{\xi_1}^2 \right)^k
           \left(
             \partial_{\eta_2} -\hf \partial_{\xi_2}^2
           \right)^{p^5-k-n}
           \left(
              \partial_{\zeta}  -\hf  \partial_{\xi_1} \partial_{\xi_2}
           \right)^n
           \nn \\
           &  \qquad \times
           \left(
               \partial_{\omega} - \frac{\xi_2}{2} \partial_{\xi_1}
               - \frac{\zeta}{2}\partial_{\eta_1} - \eta_2 \partial_{\zeta}
           \right)^{2p^5-2k-n}
    \end{align}

\end{enumerate}
\bigskip

\section{Conclusions and Outlook}
\setcounter{equation}{0}

{\color{black}

In this paper we have presented explicit expressions for the \idos\ related to the Jacobi algebra $\cg_2$.
These results are the first explicit example for $\cg_n$ with $ n \geq 2$ and elucidate that there are more variety of \idos\ than $\cg_1$ \cite{DDM}.  	
It is, therefore, natural to extend the present and that of \cite{AiDoDo} to $ \cg_n $ with $ n \geq 3. $
However, recall that the theory of parabolic Verma modules for non-semisimple Lie algebras is not developed yet.
This implies that the search for \idos\ for $\cg_n$, even for $n=2$, is a highly non-trivial problem.

The results of this paper would be useful to anyone who
would like to study the explicit implications of  $\cg_2$ invariance. Naturally, the first possible applications would be to find explicit solutions of the equations ~$\cd \cf =0$, where ~$\cd$~ would be any of the new explicit operators from Section \ref{FEidos}. Next, following the applications of the conformal and Schr\"odinger groups one may look
for explicit expressions of correlation functions invariant under $\cg_2$. Altogether, there are many possible applications.	
}

\bigskip\bigskip\bigskip 

\section*{Acknowledgments}

\nt
VKD acknowledges partial support from Bulgarian NSF Grant DN-18/1.

 \bigskip

%
%
%

\end{document}